\documentclass{amsproc}
\numberwithin{equation}{section}

\newcommand{\be}[1]{\begin{eqnarray#1}}
\newcommand{\ee}[1]{\end{eqnarray#1}} 

\newtheorem{theorem}{Theorem}[section]
\newtheorem{lemma}[theorem]{Lemma}
\newtheorem{propn}[theorem]{Proposition}      
\newtheorem{corollary}[theorem]{Corollary}  

\theoremstyle{definition}

\theoremstyle{remark}
\newtheorem{remark}[theorem]{Remark}

\newcommand{\A}{\mathcal{A}}
\newcommand{\B}{\mathcal{B}}                
\newcommand{\M}{\mathcal{M}}
\newcommand{\Ha}{\mathcal{H}}              
\newcommand{\Ru}{\mathcal{R}}
\newcommand{\U}{\mathcal{U}}                
\newcommand{\F}{\mathcal{F}}
\newcommand{\C}{\mathbb{C}}                  
\newcommand{\N}{\mathbb{N}}

\newcommand{\ad}{\mathrm{ad}}

\newcommand{\g}{\mathfrak{g}}              
\newcommand{\m}{\mathfrak{m}}

\renewcommand{\[}{{[\![}}

\newcommand{\la}{\lambda}
\newcommand{\ve}{\varepsilon}                   

\newcommand{\n}{\nonumber          }           
\newcommand{\oo}{\otimes}
\newcommand{\al}{\alpha}                         
\newcommand{\bt}{\beta}
\newcommand{\dt}{\delta}                     
\newcommand{\id}{\mbox{id}}
\newcommand{\End}{\mbox{End}}

\newcommand{\Sym}{\mbox{Sym}}
\renewcommand{\S}{\mbox{S}}              

\newcommand{\Tr}{\mbox{\rm    Tr}}     
\newcommand{\T}{\mbox{\rm     T}}
\newcommand{\si}{\sigma}
\newcommand{\gm}{\gamma}                           

\newcommand{\tr}{\triangleright}         
\newcommand{\tl}{\triangleleft}

\newcommand{\sign}{\mbox{sign}}

\begin{document}
\title[$\U_q(sl(n))$-covariant quantization of symmetric coadjoint
orbits]
{$\U_q(sl(n))$-covariant quantization of symmetric coadjoint
orbits 
via reflection equation algebra
}

\author{J. Donin}
\address{Max-Planck-Istitute f\" ur Mathematik, Vivatsgasse 7, 53111 Bonn}
\curraddr{Department of Mathematics, Bar-Ilan University, 52900 Ramat-Gan, Israel}
\email{donin@macs.biu.ac.il}
\thanks{The first author was supported in part by
Israel Academy of Sciences Grant no. 8007/99-01}

\author{A. Mudrov}
\address{Department of Mathematics, Bar-Ilan University, 52900 Ramat-Gan, Israel}
\email{mudrova@macs.biu.ac.il}
\date{August 16, 2002 and, in revised form, January 1, 2002.}
\begin{abstract}
We study relations between the two-parameter $\U_q(sl(n))$-covariant
deformation quantization on $sl^*(n)$ and the reflection equation
algebra. The latter is described by a
quantum permutation on $\End(\C^n)$ given explicitly.
The reflection equation algebra is used for
constructing the  one-parameter quantization on coadjoint
orbits, including symmetric, certain bisymmetric and nilpotent ones.
Our approach
is based on embedding the quantized function algebras on orbits 
into the algebra  of  functions  on  the  quantum
group $SL_q(n)$ by means of   reflection   equation   algebra   characters.   
\end{abstract}
\maketitle

\section{Introduction} This paper is devoted to a particular case of
the following problem.  Let $M$ be a Poisson manifold with an
action of a Lie group $G$. Let $\g$ be the  Lie algebra
of the group $G$ and $\U_h(\g)$ the corresponding quantized universal
enveloping algebra, \cite{D}. The problem is to construct a quantization 
$\A_h$ of the function algebra $\A$ on $M$ covariant with respect to the action of
$\U_h(\g)$.  In this paper, we consider the case  $\g=sl(n)$ or $gl(n)$
and $M$ is either $\g^*$ or an  orbit in $\g^*$. 
We quantize the algebra of polynomial functions on $M$ and present
it as a quotient of the corresponding tensor algebra by very natural
relations, which are similar to those in the classical case. We prove
the flatness of the deformations obtained. In that way, we quantize
all symmetric  and some bisymmetric orbits.  
Note that the analogous
relations can be written for all semisimple orbits but there arises
the question whether the quotient algebras by those relations are
flat.

Let us recall some facts about the quantization on $\g^*$ and its
orbits. It was shown in \cite{Do} that in the case $\g=sl(n)$
there exists a two-parameter deformation of the polynomial algebra
$\S(\g)$ on $\g^*$ which can be viewed as a $\U_h(\g)$-covariant
quantization of the Lie-Poisson bracket on $\g^*$. Recall that in
the classical case a natural one-parameter $\U(\g)$-covariant
quantization of $\S(\g)$ is given by the family
$\S_t(\g)=\mbox{T}(\g)[t]/J_t$, where $\mbox{T}(\g)$ is
the tensor algebra of the vector space $\g$ and the ideal
$J_t$ is generated by the elements $x\oo y - \tau(x\oo y) - t [x,y]$,
$x,y \in \g$.  Here,  $\tau$ is the flip operator on
$\g^{\oo2}$. So, the algebra $\S_t(\g)$
is quadratic-linear. By the Poincar\'e-Birkhoff-Witt
theorem, $\S_t(\g)$ is a free module over $\C[t]$.  As was proven in
\cite{Do}, \cite{Do1}, that picture can be extended to the
quantum case  for $\g=sl(n)$. Namely, there exist deformations $\tau_h$ and
$[\:\cdot,\cdot\:]_h$ of both maps $\tau$ and $[\:\cdot,\cdot\:]$ such
that the ideal $J_{h,t}$ generated by the elements $x\oo y -
\tau_h(x\oo y) - t [x,y]_h$, $x,y \in \g$, gives the two-parameter
$\U_h(\g)$-covariant quantization
$\S_{h,t}(\g)=\T(\g)[[h]][t]/J_{h,t}$
of the Lie-Poisson bracket on $\g^*$. It was also shown that
such a family does not exist for $\g\not = sl(n)$.

It is possible to prove that $\S_{h,t}(g)$ can be restricted
to any semisimple
orbit to provide a two-parameter  $\U_h(\g)$-covariant quantization
on it.
On the other hand, one can obtain a one-parameter covariant
quantization on a semisimple orbit $M$ as the subalgebra
of  quantized
functions on $G$ invariant under the action of the quantized stabilizer,
\cite{DoG}. Note that none of these statements imply
any explicit form of the quantized algebras. In this paper, we
explicitly describe
$\S_{h,t}(\g)$ and the quantum symmetric orbits in
$\g^*$. We also describe quantizations of certain nilpotent
and bisymmetric orbits. We conjecture that all bisymmetric
orbits can be quantized within our approach.

The paper is organized as follows. Sections 2 and 3 contain some
basic material about the quantum group $\U_q(\g)$, $\g=sl(n)$, and its
fundamental representation in the vector space $V=\C^n$.
In this paper, we work with the quantum group
$\U_q(\g)$ in the sense of Lusztig, \cite{Lu}, instead of $\U_h(\g)$.
The latter may be considered as the completion of $\U_q(\g)$
at the point $q=1$. Correspondingly,  instead of  $\S_{h,t}(\g)$ we
deal with the two-parameter family $\S_{q,t}(\g)$ as a free module
over $\C(q)$.
In Section~4, we study two different algebra structures on
$\T(V)\oo \T(V^*)$ covariant with respect to the actions
of either $\U^{\oo 2}_q(\g)$ or $\U_q(\g)$.  As a result, we
come to inequivalent embeddings of $\T(V\oo V^*)$ into
$\T(V)\oo \T(V^*)$.  That, in its turn, leads to different
quantizations of the polynomial algebra on $V\oo V^*$ studied in
Sections 5 and 6.  We prove that they are the reflection equation (RE)
algebra (see \cite{KSkl}) and the algebra of functions on the
quantum group known as the FRT  algebra (after Faddeev, Reshetikhin,
and Takhtajan), \cite{FRT}. It turns out that both algebras can be
described uniformly by quantum permutations on $(V\oo V^*)^{\oo2}$ but
within different categories.
In Section 7 we evaluate the lowest weight vectors of
the irreducible $\U_q(\g)$-modules in $\M^{\oo 2}$, where $\M$ is 
the algebra of $n\times n$ complex matrices.
Using this information, we derive the involutive permutation
$\tau_{RE}$ giving the commutation relations in
the reflection equation algebra. That is done in Section 8.
In Section 9, we reduce this permutation to the
submodule $\g^{\oo2}\subset \M^{\oo2}$, $\g=sl(n)$.

In Section 7, we formulate a
relation  between $\S_{q,t}(\g)$ and the RE algebra,
following \cite{Do1}. That relation involves
the reduction of $\tau_{RE}$ from $\M$ to $\g$, and
a quantum deformation of the Lie commutator
$\g^{\oo 2}\to \g$ determined by $\tau_{RE}$
up to a scalar factor.
Specifying
values of the quantum trace in the RE algebra, we obtain
one-parameter
subfamilies of $\S_{q,t}(\g)$ as quotients of the RE algebra.

Section 11 is devoted to the quantization on coadjoint
orbits in $\g^*$. We work, actually, with the adjoint
orbits in the matrix space
using the isomorphism between
adjoint and coadjoint modules and employing the relation
between the RE algebra and the two-parameter family $\S_{q,t}(\g)$.

The RE algebra appeared as an abstraction
of algebraic constructions, \cite{AFS}, \cite{KSkl}, arising from the theory
of integrable models, \cite{Cher}. It is also related
to the braid group of a solid
handlebody, \cite{K}. In this paper, we consider
the RE algebra  as the $\U_q(\g)$-covariant quantization of the polynomial
algebra on the matrix space.
It is generated by elements
arranged in the matrix $L$ and subject to a set of
quadratic relations called the reflection equation.

The conjugation transformation of $L$ with $T$, the
generating matrix of the FRT algebra whose entries
commute with the entries of $L$, is again an RE matrix,
\cite{KS}. It follows from this observation
that any character of the RE algebra specifies a homomorphism to
the FRT algebra. This is done in the same way as in
the classical geometry, where  maximal ideals of the function
algebra on $M$  correspond to points of $M$.
Each point defines an embedding of the function algebra
on the orbit into that on the group $G$.
We use this idea for constructing
quantized manifolds as quotient spaces of the
quantized group. Characters of the RE algebra
are exactly solutions to the {\em numerical} reflection
equation, \cite{KSS}.
We find such solutions in the form of projectors of rank $k<n$.
Thus, we present the one-parameter quantizations on symmetric orbits
as quotients of the RE algebra.
Simultaneously, they turns out to be subalgebras in the quantized
function algebra on the group, and this proves flatness of the 
quantizations.
There are also  solutions other than projectors,
corresponding to different paths in the parameter space of 
the two-parameter quantization $\S_{h,t}(\g)$.

We prove that, for the standard quantum linear groups,
there is an epimorphism from the RE algebra defined on
$n\times n$ matrices
onto  the RE algebra defined on $k\times k$ matrices for $k<n$. 
That enables us to build
new solutions to the matrix RE by embedding a given solution 
into a bigger matrix as the left upper block
and extending it to the whole matrix with zeros.
In particular, starting from a semisimple non-degenerate
RE matrix we gain  an additional, zero eigenvalue.
We use this method for constructing examples of quantized bisymmetric
orbits, i.e., consisting of matrices with three eigenvalues.
The complete classification of
solutions to the matrix  RE is unknown.
In particular, it is interesting to find  nilpotent matrices
providing  quantization of nilpotent orbits.
We managed to build such solutions among the matrices whose square 
is equal to zero.

\section{The quantum universal enveloping algebra $\U_q(sl(n))$}
By  $\g$ we mean the complex Lie algebra $sl(n)$, and
$\U_q(\g)$ is the corresponding quantum group. The latter is
understood in the sense of Lusztig \cite{Lu}, i.e., a free module over
the field of rational functions in $q$. In the classical limit $q\to
1$,   $\U_q(\g)$ turns into the  universal enveloping algebra
$\U(\g)$.  The quantum group $\U_q(\g)$ is  generated by the elements
$H_i$, $X^\pm_i$, $i=1,\ldots,n-1,$ satisfying the commutation
relations
$$
[H_i,X^\pm_i]=\pm 2 X^\pm_i,\quad [H_i,X^\pm_{i\pm 1}]=
               \mp X^\pm_{i\pm 1},\quad
$$
$$
[H_i,X^\pm_j]=0, \quad |i-j|>1,
$$
$$
[X^+_i,X_j^-]=\delta_{ij}\frac{q^{H_i}-q^{-H_i}}{q-q^{-1}}.
$$
Besides, the Serre relations hold:
$$
(X_i^\kappa)^2 X^\kappa_{i\pm 1}- 2_q X^\kappa_i X^\kappa_{i\pm 1} X^\kappa_i +
X^\kappa_{i\pm 1} (X^\kappa_i)^2 =0,
\quad [X^\kappa_i, X^\kappa_j]=0, \;|i-j|>1,
$$
where $\kappa=\pm$.
Here $2_q = q+q^{-1}$. In general, the quantum integer numbers
are defined as $n_q= \frac{q^n-q^{-n}}{q-q^{-1}}$.
The coproduct $\Delta$, counit $\ve$, and antipode $\gm$ are
$$
\Delta(H_i)=H_i\oo 1 + 1\oo H_i, \quad
\Delta(X^\pm_i)=X^\pm_i\oo q^{-\frac{H_i}{2}} + q^{\frac{H_i}{2}}\oo X^\pm_i.
$$
$$
\ve(H_i)=\ve(X^\pm_i)=0,\quad \gamma(H_i)=-H_i,\quad \gamma(X^\pm_i)= - q^{\mp 1}X^\pm_i.
$$
In the (completed) tensor square of $\U_q(\g)$, there is an element $\Ru$ called
the universal R-matrix. It satisfies the conditions
\be{}
(\Delta\oo \id)(\Ru)=\Ru_{13}\Ru_{23}
, & &
(\id\oo \Delta)(\Ru)=\Ru_{13}\Ru_{12},
\label{R1}
\\
\Ru \Delta(x) &=& \Delta'(x)\Ru, \quad x \in \U_q(\g)
\label{R2}
.
\ee{}
The subscripts distinguish the legs of tensor objects and the prime
stands for the opposite coproduct.
We adopt the Sweedler symbolic notations $\Delta(x)=x_{(1)}\oo x_{(2)}$
for the coproduct
and  denote the antipode by the
bar: $\gamma(x)=\bar x $,  $x\in \U_q(\g)$.
The exact expression for  $\Ru$ can be found in \cite{Ros,KR,KT}.
We do not use it in the sequel, rather the image of $\Ru$ in the
fundamental representation on $\C^n$.
\section{Fundamental representation of $\U_q(sl(n))$}
\label{SecFR}
By $\M$ we denote the algebra of $n\times n$ complex matrices with the
multiplication $e^i_j e^k_l =  e^i_l\delta^k_j$ defined on the matrix
units $e^i_j$.
Here $\delta^k_j$ stand for  the Kronekker symbol.
The elements of $\M$  are considered as right
endomorphisms of $V=\C^n$ with the action $e_i e^k_l =
e_l\delta^k_i$
on the  basis elements $e_i\in V$.
Formulas
$$
\rho(H_i) = e^i_i - e^{i+1}_{i+1}, \quad
\rho(X^+_i) = e^{i}_{i+1}        , \quad
\rho(X^-_i) = e^{i+1}_{i}, \quad i=1,\ldots, n-1.
$$
define a homomorphism $\rho$ of the algebras $\U_q(\g)$ into $\M$.
On the Chevalley generators, it is given by the same formulas
as in the classical limit $q\to 1$.

The homomorphism $\rho$ defines the right {\em adjoint}
and left {\em coadjoint} actions of  $\U_q(\g)$ on $\M$:
\be{}
\label{ad}
A \tl \ad_\rho(x)  &= &\rho(\overline {x_{(1)}})A\rho({x_{(2)}})
,
\\
\label{coad}
\ad^*_\rho(x) \tr A& = & \rho({x_{(2)}})A\rho(\overline {x_{(1)}})
\ee{}
for $x \in \U_q(\g)$ and $A\in \M$. These actions are conjugate via
the trace pairing on $\M$:
$$
\Tr\biggl(\bigl(A \tl \ad_\rho(x)\bigr)\; B\biggr)
=\Tr\biggl(A \;\bigl(\ad^*_\rho(x) \tr B\bigr)\biggr).
$$
The left coadjoint module $\M$ equipped with action
(\ref{coad}) will be  denoted $\M^*$.

The element
\be{}
\label{R}
R=
q^{-\frac{1}{n}}(q \sum_i e^i_i\oo e^i_i + \sum_{i\not = j} e^i_i\oo e^j_j +
\omega \sum_{i<k}e^k_i \oo e^i_k) \in \M\oo \M,
\ee{}
where $\omega= q-q^{-1}$,
satisfies the Yang-Baxter equation
$$
R_{12}R_{13}R_{23}=R_{23}R_{13}R_{12}.
$$
It differs from the R-matrix used in \cite{FRT} by the overall
factor $q^{-\frac{1}{n}}$. In this form, this is the image of
the universal R-matrix $\Ru$ under the homomorphism $\rho^{\oo2}$.
Let $\si_{V\oo V}$ be the ordinary (classical) permutation on $V\oo V$.
The braid matrix $S=q^{\frac{1}{n}}\si_{V\oo V} R$ 
satisfies the Hecke condition $S^2=\omega S+1$ and is
represented by the sum of two orthogonal projectors 
\be{}
P^+_q=\frac{1}{2_q}(q - S)
,&\quad&
P^-_q=\frac{1}{2_q}(q^{- 1}+ S),
\label{Ppm}
\ee{}
The matrix $S$ commutes with all elements $(\rho\oo \rho)\Delta(x)$,
$x\in \U_q(\g)$.


Let $V^*$ be the space of linear functionals on $V$. We choose a
basis $\{e_i\}$ in $V$ and denote $\{f^i\}$ its dual
with respect to the canonical pairing between $V$ and $V^*$.
\begin{propn}
\label{bacis}
The homomorphism $\rho\colon \U_q(\g)\to \M$ induces
the right and left actions of $\U_q(\g)$ on
the tensor algebra $\T(V\oplus V^*)$.
They are given on the basis elements $e_i\in V$ and $f^i\in V^*$
by
\be{}
e_i \tl x  = \sum_{\al=1}^n e_\al \rho(x)^\al_i 
, &\quad &
f^i \tl x = \sum_{\al=1}^n f^\al \rho(\bar x)_\al^i, 
\label{ractions} \\
x\tr e_i = \sum_{\al=1}^n e_\al \rho(\bar x)^\al_i ,& \quad & x\tr f^i =
\sum_{\al=1}^n f^\al  \rho(x)_\al^i,  
\label{actions}  
\ee{}  
for  $x\in
\U_q(\g)$. Formula  $\langle  e_i,  f^j\rangle=\delta^j_i$   defines   a
$\U_q(\g)$-invariant right pairing between  $V$  and  $V^*$.  The  right
module $V^*\oo  V$  with  the  action  induced  by  (\ref{ractions})  is
isomorphic to the adjoint module $\M$ via the correspondence $f^i\oo e_j
\to e^i_j$. The left module  $V\oo  V^*$  with  the  action  induced  by
(\ref{actions}) is isomorphic to the coadjoint  module  $\M^*$  via  the
correspondence $e_j \oo f^i \to e^i_j$. 
\end{propn} 
\begin{proof} 
Direct verification. 
\end{proof}

One  can  introduce  the   following   $\U_q(\g)$-invariant   involutive
permutations among the elements of $V$ and $V^*$. 
\begin{propn} 
The involutive permutations
\be{}
\tau_{V\oo  V}(e_i\oo  e_j)
&=&
\left\{
\begin{array}{l}
\sign(j-i)\frac{1-q^2}{1+q^2} e_i\oo  e_j  +\frac{2q}{1+q^2} e_j\oo  e_j
,\quad  i\not =j,\\
e_i\oo e_i
,\;  i =j, 
\end{array}
\right.
\label{VV}   
\\
\tau_{V^*\oo V^*}(f^i\oo  f^j)
&=&
\left\{
\begin{array}{l}
\sign(i-j)\frac{1-q^2}{1+q^2} f^i\oo  f^j  +\frac{2q}{1+q^2} f^j\oo  f^j
,\;  i\not =j,\\
f^i\oo f^i
,\quad  i =j, 
\end{array}
\right.
\label{V*V*}   
\\
\tau_{V\oo  V^*}  (e_j\oo  f^i)  
&=&
\left\{
\begin{array}{l}
q^{-\frac{1}{n}}
(q  f^i\oo  e_i  +\omega
\sum_{k>i}f^k\oo e_k)
,\;  i =j,\\
q^{-\frac{1}{n}}
f^i\oo e_j
,\quad  i\not =j, 
\end{array}
\right.
\label{VV*}   
\ee{}
are invariant with respect to the left $\U_q(\g)$-action (\ref{actions}).
The permutation 
\be{}
\tau^r_{V\oo  V}(e_i\oo  e_j)
&=&
\left\{
\begin{array}{l}
\sign(i-j)\frac{1-q^2}{1+q^2} e_i\oo  e_j  +\frac{2q}{1+q^2} e_j\oo  e_j
,\quad  i\not =j,\\
e_i\oo e_i
,\quad  i =j, 
\end{array}
\right.
\label{rVV}     
\ee{}
is invariant with respect to the right $\U_q(\g)$-action (\ref{ractions}).
\end{propn} 
\begin{proof} 
As a lenear operator, composition of $\Ru$ acting on $V\oo V$ and $V^*\oo V^*$
with the ordinary flip has two eigenspaces.
They are deformations
of symmetric and antisymmetric 2-tensors, and involutions  (\ref{VV}), 
(\ref{V*V*}), and  
(\ref{rVV})  are defined as multiplication by $\pm 1$ on those subspaces.
Operation (\ref{VV*})  is  readily
obtained from the action of the universal R-matrix:
$\tau_{V^*\oo  V}(f\oo e) = \Ru_2\tr e \oo \Ru_1\tr f$,
$\tau_{V\oo  V^*}(e\oo f) = \Ru^{-1}_1\tr f \oo \Ru^{-1}_2\tr e$,
$e\in V$, $f\in V^*$. To compute  $\tau_{V\oo  V^*}$,  it is enough
to know only the image (\ref{R}) of $\Ru$ in $\M^{\oo 2}$. 
\end{proof}  

\section{Algebra  $\T(V)\oo\T(V^*)$  in  braided categories} 
\label{REA-P} 
In  this  section,  we  consider  the  tensor
algebras $\T(V)$ and $\T(V^*)$ from the different points  of  view:  as  
objects  from  the   category   of   either   $\U_q(\g)$-   or   $\U^{\oo
2}_q(\g)$-representations. Correspondingly,  the  algebra  structure  on
$\T(V)\oo\T(V^*)$ may be introduced  in  different ways, as the tensor product
of algebras in the those categories. A
particular choice of  the  category  determines  embedding  of  
$\T(V\oo V^*)$,  the   tensor   algebra   of   the   space   $V\oo   V^*$,   
into $\T(V)\oo\T(V^*)$. That leads to different  quantizations  of  
the polynomial algebra on $V\oo V^*$.

Let us remind the construction of the tensor product of algebras in the category
of a quantum group modules. Let $\Ha$ be a quantum group with
an R-matrix $\Ru$.
Recall that a (left) $\Ha$-module
algebra $\A$ is an associative
algebra with unit in the category of $\Ha$-representations:
\be{}
\label{ma}
x\tr (ab) & = & (x_{(1)}\tr a) (x_{(2)}\tr b),
\quad x\tr 1 = \ve(x) 1,
\quad x \in \U_q(\g),
\; a,b  \in \A.
\ee{}
Given two module algebras $\A$ and $\B$, one can
introduce their (braided) tensor product $\A\check\oo \B$, which is
again an $\Ha$-module algebra.
\begin{propn}
Let $\Ru_1 \oo \Ru_2$ denote the decomposition
of the universal R-matrix into the two tensor factors (summation supressed).
The formula
\be{}
\label{tp}
(a_1\check\oo b_1)(a_2\check\oo b_2)=
a_1(\Ru_2\tr a_2)\check\oo (\Ru_1\tr b_1)b_2,\quad a_i\in \A, 
\quad b_i \in \B
\ee{}
defines an associative multiplication on
$\A\check\oo \B$ turning it into
an $\Ha$-module algebra.
Embeddings $\A\to \A \check\oo 1$ and $\B\to 1\check \oo \B$ are
homomorphisms of algebras.
\end{propn}
\begin{proof}
Associativity of multiplication (\ref{tp}) follows from
(\ref{R1}). Compatibility with the action of $\Ha$ is a
consequence of  (\ref{R2}). The last statement of the proposition is
immediate due to the equality $(\ve\oo \id)(\Ru)=  (\id\oo
\ve)(\Ru)= 1\oo 1$ following from the definition of the universal
R-matrix.
\end{proof}
Note that multiplication (\ref{tp}) is characterized by the
permutation relation
\be{} ba   =  (\Ru_2\tr a)(\Ru_1\tr b)  \label{ab=Rba}
\ee{}
between the elements $a\in\A$ and $b\in \B$.

\begin{propn}
\label{TVV*} 
Let $V$ and $U$ be $\Ha$-modules.
The subalgebra in $\T(V)\check\oo\T(U)$ generated
by the submodule $V\check \oo U$ is isomorphic to
the tensor algebra  $\T(V\oo U)$.
It is an algebra in the category of
$\Ha$-modules.
\end{propn}
\begin{proof}
Embedding $V\oo U \to \T(V)\check\oo\T(U)$
is extended
to a homomorphism of the free algebra $\T(V\oo U)$.
It is an monomorphism due  to invertibility of permutation
(\ref{ab=Rba}).
The subalgebra generated by $V\oo U$ is $\Ha$-invariant,
so it is a module algebra.
\end{proof}
\noindent
We will use this proposition when $U=V^*$.

Let us  apply the construction above to the algebras $\A=\T(V)$,  $\B=\T(V^*)$
considered as $\U_q(\g)\oo \U_q(\g)$- and $\U_q(\g)_{op}\oo \U_q(\g)$-module algebras in 
the following way.
First let us note that a right action of a Hopf algebra is the same as
a left one for its opposite.
Right and left actions (\ref{ractions}) and (\ref{actions}) of the algebra 
$\U_q(\g)$  are extended to the left actions 
of the Hopf algebras $\U_q(\g)\oo \U_q(\g)$- and $\U_q(\g)_{op}\oo \U_q(\g)$
on $V$ and $ V^*$:
\be{}
\label{UU}
(x\oo y) \tr v = \ve(y) x\tr v, &\quad &
(x\oo y) \tr \phi = \ve(x) y\tr \phi,
\ee{}
where $v \in V$ and $\phi \in V^*$. 
These actions  are extended to the actions on
$\T(V)$, $\T(V^*)$, and on the tensor product  $\T(V)\oo\T(V^*)$.
\begin{propn}
\label{UUtp}
The algebra $\T(V)\check\oo \T(V^*)$ in the categories of
$\U_q(\g)\oo \U_q(\g)$- and $\U_q(\g)_{op}\oo \U_q(\g)$-modules has 
the multiplication:
\be{}
(a_1\check\oo b_1)(a_2\check\oo b_2)=
a_1 a_2\check\oo  b_1b_2,\quad a_i\in \T(V), \quad b_i \in  \T(V^*).
\ee{}
\end{propn}
\begin{proof}
In both cases, the universal R-matrix is the product
of two $\Ru$-s supported in the two different copies
of  $\U_q(\g)$. By construction (\ref{UU}), each copy acts
trivially either on
$\T(V)$ or on $\T(V^*)$, so the universal R-matrix turns into
the identity operation in (\ref{tp}).
\end{proof}

Note that, being considered in the category of $\U_q(\g)$-modules,
the tensor product $\T(V)\check\oo \T(V^*)$ has a non-classical multiplication
characterized by the permutation relation (\ref{ab=Rba}), where $\Ru$ is
the R-matrix for $\U_q(\g)$.

Due to Proposition \ref{TVV*}, the tensor algebra $\T(V\oo V^*)$
has different realizations in the categories
of either $\U_q(\g)_{op}\oo\U_q(\g)$- or $\U_q(\g)$-representations.
One can introduce additional relations in $\T(V\oo V^*)$ resulting
in two different quantizations of the polynomial algebra on
$V\oo V^*$; they are the FRT and RE algebras. We consider them in
Sections \ref{SecFRT} and \ref{SecRE}  and present their defining
commutation relations in terms of involutive permutations $\tau_{FRT}$
and $\tau_{RE}$.

Taking into account Proposition \ref{UUtp}, we
reserve the symbol $\check \oo$ only for the tensor product in  the
category of $\U_q(\g)$-modules where the permutation (\ref{ab=Rba})
between $V$ and $V^*$ is non-classical. So, by $\T(V)\check\oo \T(V^*)$
we will denote the tensor product considered as
$\U_q(\g)$-module.

\section{Algebra $\A_{FRT}(\M)$} \label{SecFRT}
Denote $\A_{FRT}(\M)$ be the associative
unital algebra over $\C(q)$ generated by the matrix elements $T^i_j$
subject to relations \be{} S T_1 T_2 & = & T_1 T_2 S.
\label{RTT=TTR}
\ee{}
Here, $S$ is the
Hecke matrix defined in Section \ref{SecFR}.
The algebra $\A_{FRT}(\M)$ was introduced in \cite{FRT} as a quantized
polynomial algebra on the space of matrices. It is endowed with
the $\U_q(\g)_{op}\oo \U_q(\g)$-module structure coming from that on $V\oo
V^*$.  Relations (\ref{RTT=TTR}) are given by a quantum permutation
$\tau_{FRT}$ on the space $V\oo V^*$, which we are going to present
explicitly.

Consider $\T(V\oo V^*)$ as a subalgebra of $\mbox{T}(V)\oo
\mbox{T}(V^*)$ in the tensor category of $\U_q(\g)_{op}\oo \U_q(\g)$-modules,
along the line of Proposition \ref{TVV*}.  
Consider the permutation
\be{}
\label{tauFRT}
\tau_{FRT}&=&\si_{V\oo V^*}\circ\tau^r_{V\oo V}\circ
\tau_{V^*\oo V^*}\circ \si^{-1}_{V\oo V^*}
\ee{}
on $(V\oo V^*)^{\oo 2}$, which is the composition of elementary 
permutations (\ref{V*V*}),
(\ref{rVV}), and the classical flip $\si_{V\oo V^*}$. The permutation 
$\tau_{FRT}$ is $\U_q(\g)_{op}\oo \U_q(\g)$-invariant.
Let  $(V\oo V^*)^{\oo2}=I_{FRT}^+\oplus I_{FRT}^-$ be the decomposition 
into symmetric and skew-symmetric submodules with respect to $\tau_{FRT}$.
It is easy to see that $I_{FRT}^-$ is a deformation of the exterior square
$(V\oo V^*)^{\wedge2}$.  

\begin{theorem} \label{FRT} Let ${\mathcal J}_{FRT}$ be  the ideal in  $\T(V\oo V^*)$ generated by
$I_{FRT}^-$. The quotient algebra $\T(V\oo V^*)/{\mathcal J}_{FRT}$
is isomorphic to $\A_{FRT}(\M)$.  
\end{theorem}
\begin{proof}
Condition (\ref{RTT=TTR}) is equivalent to the pair of
relations $P^\pm_q T_1T_2 P^\mp_q=0$, where $P^\pm_q$ are the 
projectors (\ref {Ppm}).
Let us take a basis $\{e_i\}\subset V$ and  its right dual
$\{f^i \}\subset V^*$,
i.e.,
$\langle e_i, f^j \rangle=\delta^j_i$.
The elements  $e_i\oo f^j$ transform as entries of the matrix $T$
under the action of $\U_q(\g)_{op}\oo \U_q(\g)$, so we can identify them.  The
product of two matrix entries  can be rewritten as $T^i_j T^m_n =
(e_j\oo f^i)(e_n\oo f^m) = e_j e_n \oo f^i f^m$, since the elements of
$V$ and $V^*$ commute.  This implies that the quadratic submodule
$I_{FRT}^-$ goes over into the submodule spanned by
$\sum_{\al,\bt,\mu\nu}(P^\pm_q)_{jn}^{\al\bt} e_\al e_\bt \oo 
f^\mu f^\nu (P^\mp_q)_{\mu\nu}^{im}$,  under the embedding 
$T(V\oo V^*) \to T(V)\oo T(V^*)$.
This submodule is formed by symmetric tensors with respect to 
the permutation $\tau^r_{V\oo V}\circ \tau_{V^*\oo V^*}$.
\end{proof}
\section{Algebra $\A_{RE}(\M)$}
\label{SecRE}
The reflection equation algebra $\A_{RE}(\M)$
is  generated by the entries of the matrix $L$ subject to the
quadratic relations, \cite{AFS, KSkl},
\be{} \label{RER}
L_2 S L_2 S &=& S L_2 S L_2.
\ee{}
Here, $S$ is the  Hecke matrix introduced in Section \ref{SecFR}.
The  generators $L^i_j$
form the coadjoint module $\M^*$ via the correspondence
$L^i_j\to e^i_j$. Action (\ref{coad}) is extended over
$\A_{RE}(\M)$ turning it into a $\U_q(\g)$-module algebra;

The element $Q=\Ru_{21}\Ru$ defines a morphism of $\U_q(\g)$-modules  
$\A_{RE}(\M)\to\U_q(\g)$ by the correspondence 
$L\to Q_2\oo \rho(Q_1)$, \cite{Mj},
where $\U_q(\g)$ is considered as a left $\U_q(\g)$-module 
via the adjoint action.
The  algebra $\A_{RE}(\M)$ is a flat deformation of the polynomial
algebra on the matrix space $\M$, \cite{Mj}.  

The  algebra $\A_{RE}(\M)$ can be described 
in the same way as $\A_{FRT}(\M)$ by a
quantum permutation $\tau_{RE}$ .
Following Proposition \ref{TVV*}, consider the tensor algebra $\T(V\oo
V^*)$ as the subalgebra $\T(V\check\oo V^*)\subset \T(V)\check \oo
\T(V^*)$ in the category of $\U_q(\g)$-modules.  
Introduce the involutive operation 
\be{}
\tau_{RE}&=&\tau_{V\oo V^*}\circ\tau_{V\oo V}\circ 
\tau_{V^*\oo V^*}\circ \tau^{-1}_{V\oo V^*},
\ee{}
which is the composition of
elementary permutations (\ref{VV}--\ref{VV*}).
The permutation $\tau_{RE}$ is $\U_q(\g)$-invariant.
Consider the decomposition 
$(V\oo V^*)^{\oo2}=I_{RE}^+\oplus I_{RE}^-$ 
into the symmetric and skew-symmetric submodules with respect to $\tau_{RE}$.
It is easy to see that the submodule $I_{RE}^-$ is a deformation of the 
exterior square $(V\oo V^*)^{\wedge2}$.

\begin{theorem} \label{REA} Let ${\mathcal J}_{RE}$ be  the ideal 
in  $\T(V\oo V^*)$ generated by $I_{RE}^-$.
The quotient algebra
$\T(V\check\oo V^*)/{\mathcal J}_{RE}$ is isomorphic to  $\A_{RE}(\M)$.
\end{theorem}
\begin{proof}
As in the proof of Theorem \ref{FRT}, let $\{f^i\}\subset V^*$ be
the right dual to the basis $\{e_i\}\subset V$.
Due to  Proposition \ref{bacis}, the elements
$e_i\check \oo f^j$
form the coadjoint module $\M^*$ for $U_q(\g)$, so we can put
$L^j_i=e_i \check \oo f^j$.  By Proposition \ref{TVV*}, the algebra
$\T(V\check\oo V^*)$ is isomorphic to the free associative algebra
generated by $L^i_j$.  It remains to show that the quadratic submodule
generating the ideal ${\mathcal J}_{RE}$ is isomorphic to the submodule specified by
the reflection equation. They are conjugate by the
$\tau_{V\oo V^*}$, which is the restriction to $V\oo V^*$ of  permutation
(\ref{ab=Rba}).
Consider the equality
$e_n e_i \check \oo f^k f^n \oo e^n_k \oo e^i_m = \si_{V\oo V} L_2 S L_2$ 
in the algebra $T(V)\check \oo T(V^*)$.
Since $T(V)$ is considered as a left module for $\U_q(\g)$, the 
permutation $\tau_{V\oo V}$ is $\si_{V\oo V}$-conjugate with the 
permutation  $\tau^r_{V\oo V}$.
Therefore, the quadratic submodule of q-symmetric tensors 
in $V^{\oo2} \oo V^{*\oo2}$ is $\si_{V\oo V}$-conjugate with
the $\U_q(\g)_{op}\oo\U_q(\g)$-invariant submodule
providing the algebra $\A_{FRT}(\M)$:
$\si_{V\oo V}S \si_{V\oo V} \bigl(\si_{V\oo V} L_2 S L_2\bigr) =
\bigl(\si_{V\oo V} L_2 S L_2\bigr)S$.
This is equivalent to the reflection equation (\ref{RER}).
\end{proof}
\begin{remark} In the classical case, the
polynomial algebra on $V\oo V^*$ can be built equivalently
as a $\U(\g)$- or $\U^{\oo 2}(\g)$-module. In the quantum case, depending
on the point of view,
one comes to  either  FRT or RE algebras. 
It can be shown that there is a sequence of twist transformations
in the quasitensor category of $\U_q(\g)\oo \U_q(\g)$-modules
relating algebras $\A_{FRT}(\M)$ and $\A_{RE}(\M)$.
This is a subject of our forthcoming publication.
\end{remark}

\section{Irreducible submodules of $\M^{\oo2}$}
\label{weights}
This technical section contains information about
submodules of $\M^{\oo2}$ and its dual in terms of highest
and lowest weight vectors. They will be used in computing the
quantum permutation $\tau_{RE}$ on $(V\oo V^*)^{\oo 2} \sim \M^{*\oo2}$.
It is convenient to evaluate the conjugate operator
$\tau^*_{RE}$ in the dual space  first.
For that reason, we  consider
three $\U_q(\g)$-module structures on  $\M^{\oo 2}$:
\be{}
(A\oo B)\tl x
& = &\rho(\overline{x_{(2)}})A \rho(x_{(3)}) \oo
\rho(\overline{x_{(1)}})B \rho(x_{(4)})
\label{al},\\
(A\oo B)\tl x
& = &\rho(\overline{x_{(1)}})A \rho(x_{(2)}) \oo
\rho(\overline{x_{(3)}})B \rho(x_{(4)})
\label{bt},\\
x \tr (A\oo B)
& = &\rho(x_{(2)})A \rho(\overline{x_{(1)}}) \oo
\rho(x_{(4)})B \rho(\overline{x_{(3)}})
\label{dt}.
\ee{}
Actions (\ref{al}) and (\ref{bt}) are nothing else than
$\ad_{\rho\oo \rho}$ and $\ad_\rho \oo \ad_\rho$, correspondingly.
The left action (\ref{dt}) is conjugate to (\ref{bt}) with respect to  the
ordinary trace pairing and it is just $\ad^*_\rho\oo \ad^*_\rho$.

A basis of the lowest weight vectors of irreducible $\U_q(\g)$-submodules with respect
to  the action (\ref{dt}) may be calculated directly. It is:

\be{}
\dt_1&=& \sum_{i=1}^n q^{-2i+1}e^i_i\oo e^n_1 ,    
\n\\
\dt_2 &=& \sum_{i=1}^n e^n_1\oo q^{-2i+1}e^i_i ,
\n\\
\dt_3 &=&   \sum_{i=1}^n  e^i_1\oo e^n_i ,
\n\\
\dt_4 &=&  \sum_{i=1}^n a_i e_i^n\oo e_1^i,
\quad
a_i=1,q^{-3},\ldots,q^{-2n+3},q^{-2n},
\n\\
\dt_5 &=& \sum_{i,k=1}^n  q^{-(2i-1)}q^{-(2k-1)}e^i_i\oo e^k_k ,
\n\\
\dt_6 &=& \sum_{i,k=1}^n  a_ie^k_i\oo e^i_k,
\quad
a_i=q^{-(2n-1)},q^{-(2n-3)},\ldots,q^{-3},q^{-1},
\n\\
\dt_7 &=&
qe_1^n\oo e_2^{n-1} + q^{-1} e_2^{n-1} \oo e_1^n
-  e_1^{n-1}\oo e_2^n -          e_2^n \oo e_1^{n-1},
\n\\
\dt_8 &=&
e_1^n\oo e_1^n,
\n\\
\dt_9 &=&
e_1^n\oo e_1^{n-1} - q^{-1} e_1^{n-1} \oo e_1^n ,
\n\\
\dt_{10} &=&
e_1^n\oo e_2^n - q^{-1} e_2^n \oo e_1^n,
\n
\ee{}
The coefficients $a_i$ entering $\dt_4$ and $\dt_6$ decrease by $q^{-2}$ each step
within the interval $i=1,\ldots, n$.
Here, we assume $n>2$. The case $n=2$ is simple and left to the reader.

The vectors $\dt_5$ and $\dt_6$ are invariant elements under action (\ref{dt}). They are
symmetric with respect to the permutation $\tau^*_{RE}$, since 
that is so for the invariant elements
in the classical situation.
The vectors  $\dt_7$ and $\dt_8$ belong to the
symmetric part, too, while $\dt_9$ and  $\dt_{10}$ are antisymmetric with respect to
$\tau^*_{RE}$. Indeed, 
they generate submodules of multiplicity one in $\M^{\oo2}$ and are deformations
of classical submodules belonging to symmetric and antisymmetric parts, respectively.
The situations is more
complicated in what concerns vectors $\dt_1,\ldots,\dt_4$ because
they generate the isotypic component of type $\g$ and, in the
limit $q=1$, the corresponding classical component 
intersects the  
symmetric and skew-symmetric parts of $\M^{\oo2}$ with multiplicity two.  
So the problem 
of evaluating the permutation $\tau_{RE}$ reduces
to calculating its restriction to the $\g$-isotypic component of
$\M^{*\oo2}$.

First,
we compute the dual conjugate operation to $\tau_{RE}$ on
the right module $\M^{\oo 2}$ with action (\ref{bt}).  A basis of
lowest weight vectors has the dual consisting of highest weight ones;
by that reason we also consider the highest weight vectors  $\al^i$
and $\bt^i$ of type $\g$ with respect to  actions (\ref{al}) and
(\ref{bt}), respectively:
\be{}
\al^1 &=&  1\oo e^1_n , \n\\ \al^2
&=&  \sum_{i=1}^n a_i e^1_n\oo e^i_i, \quad a_i=q,1,\ldots,1,q^{-1} ,
\n\\
\al^3 &=&  \sum_{i=1}^n a_i e^1_i\oo e^i_n,  \quad a_i=q,1,\ldots,1,q^{-1},
\n\\
\al^4 &=&  \sum_{i=1}^n a_i e^i_n\oo e^1_i,  \quad a_i=1,\ldots,1,q^{-1},
\n\\
\bt^1 &=&  1\oo e^1_n ,
\n\\
\bt^2 &=&  e^1_n\oo 1 ,
\n\\
\bt^3 &=&  \sum_{i=1}^n a_i e^1_i\oo e^i_n,
\quad
a_i=q^{-(2n-1)},q^{-(2n-3)},\ldots,q^{-3},q^{-1},
\n\\
\bt^4 &=&  \sum_{i=1}^n a_i  e^i_n\oo e^1_i, \quad a_i=q,1,\ldots, 1,q^{-1}.
\ee{}
Let us denote the isotypic $\g$-type components with
respect to actions (\ref{al}), (\ref{bt}), and (\ref{dt}) by
$\M^{\oo2}_\al$, $\M^{\oo2}_\bt$, and  $\M^{*\oo2}_\dt$. They are
generated by the vectors $\al^i$, $\bt^i$, and $\dt_i$, $i=1,\ldots,
4$, respectively.

\section{The quantum permutation $\tau_{RE}$}
\label{qsl}
In this section, we calculate the permutation $\tau_{RE}$
in terms of the lowest weight vectors introduced in
the previous section.
The permutations $\tau_{FRT}$ and $\tau_{RE}$ define quadratic relations
in the deformed algebras of functions
on matrices.
In both cases, relations (\ref{RTT=TTR}) and (\ref{RER})
are formulated using the algebra
structure on $\M^{\oo2}$, i.e., in the dual setting.
Thus, the FRT and RE algebra  permutations
are introduced through the dual operations $\tau^*_{FRT}$ and
$\tau^*_{RE}$ conjugate to  $\tau_{FRT}$ and $\tau_{RE}$
by the trace pairing.
The involutions $\tau^*_{FRT}$ and $\tau^*_{RE}$
are invariant with respect to right actions
(\ref{al}) and (\ref{bt}).
It is natural to compute them first,
using the algebra structure on $\M^{\oo2}$,
and then evaluate $\tau_{RE}$ by duality.

The subspace of  $\Omega \in \M^{\oo2}$ such that $S\Omega=\Omega S$
is the annulator of the submodule defining  FRT relations
(\ref{RTT=TTR}).  Therefore, the involution $\tau^*_{FRT}$ is
determined by the conjugation transformation with the Hecke matrix
$S$.
Solutions to the equation $S\si(\Omega)=\si(\Omega)S$, where the
map $\si$ is defined as $ \Omega \to (1\oo \Omega_1) S (1\oo
\Omega_2)$, $\Omega=\Omega_1\oo\Omega_2\in \M^{\oo2}$, form the
submodule annihilating the RE relations.
Hence, the involutions $\tau^*_{RE}$ and
$\tau^*_{FRT}$ are $\si$-conjugate:
$\tau^*_{RE}=\si^{-1}\tau^*_{FRT}\si$.
Since the Hecke matrix $S$ is invariant under the adjoint action
$\ad_{\rho\oo\rho}$, the map $\si$ intertwines actions
(\ref{bt}) and (\ref{al}).  It is an isomorphism of
modules, being a deformation
of the classical flip  $\si_{\M\oo \M}$.

As was noted in Section \ref{weights},
only the restriction of $\tau^*_{RE}$ to the isotypic $\g$-type
component requires a special consideration. 
Restricted to the $\g$-component $\M^{\oo2}_\al$, the conjugation
operation with the matrix $S$ has eigenvalues $+1$ and $-q^{\pm 2}$.
The eigenspace $\M^{\oo2}_{\al;+}$ corresponding to
the eigenvalue $+1$ consists of symmetric tensors, in the sense
of the permutation $\tau^*_{FRT}$. Eigenvectors of the
eigenvalues $-q^{\pm 2}$ may be called
$\tau^*_{FRT}$-antisymmetric.
We denote this submodule $\M^{\oo2}_{\al;-}$.
Thus, the conjugation with the matrix $S$ determines decomposition
\be{}
\label{alDec}
\M^{\oo2}_\al=\M^{\oo2}_{\al;+}\oplus\M^{\oo2}_{\al;-},
\ee{}
where both $\M^{\oo2}_{\al;\pm}$ contain irreducible
submodule $\g$ with multiplicity two.
Decomposition (\ref{alDec}) induces decompositions
\be{}
\label{btDec}
\M^{\oo2}_\bt=\M^{\oo2}_{\bt;+}\oplus\M^{\oo2}_{\bt;-},
\\
\label{siDec}
\M^{*\oo2}_\dt=\M^{*\oo2}_{\dt;+}\oplus\M^{*\oo2}_{\dt;-}.
\ee{}
into the symmetric and antisymmetric parts.
The subspaces $\M^{\oo2}_{\bt;\pm}$ are the images of
$\M^{\oo2}_{\al;\pm}$ via the inverse transformation $\si^{-1}$.
They yield the involution $\tau^*_{RE}$.
The subspaces $\M^{*\oo2}_{\dt;\pm}$ are the annulators
of $\M^{\oo2}_{\bt;\mp}$ with respect to the trace pairing
between  $\M^{\oo 2}$ and  $\M^{*\oo 2}$. The permutation
$\tau_{RE}$ is determined by the submodules
$\M^{*\oo2}_{\dt;\pm}$.
\begin{propn}
\label{tauRE}
The $\g$-type submodules of $\tau_{RE}$-symmetric and antisymmetric tensors
are generated by the following highest and lowest weight vectors:
$$
\begin{array}{ll}
\M^{\oo2}_{\al;+} \; : \;\al^3+\al^4, \quad \al^1+\al^2+\omega \al^3
,& \;
\M^{\oo2}_{\al;-} \; : \; q^{\pm 1}\al^1 - q^{\mp 1}\al^2 \mp \al^3 \pm \al^4
,\\[12pt]
\M^{\oo2}_{\bt;+} \; : \; \bt^1+\bt^2-\omega\bt^3,\quad \bt^3+\bt^4+\omega\bt^1
,& \;
\M^{\oo2}_{\bt;-} \; : \; \bt^1-\bt^2            ,\quad   \bt^3-\bt^4
,\\[12pt]
\M^{*\oo2}_{\dt;+} \;: \; \dt_3+\dt_4-\omega\dt_2,\quad  \dt_1+\dt_2+\omega\dt_3
,& \;
\M^{*\oo2}_{\dt;-} \;: \; \dt_1-\dt_2            ,\quad   \dt_3-\dt_4
.
\end{array}
$$
\end{propn}
\begin{proof}

The eigenvectors of the conjugation by $S$ restricted to
$\M^{\oo2}_\al$
are computed directly. Thus we obtain  (\ref{alDec}).
The transformation $\si$ acts  on the  highest weight vectors of
the $\g$-type colmponents as
$$
\si\colon \left(
\begin{array}{c}
\bt^1 \\
\bt^2 \\
\bt^3 \\
\bt^4
\end{array}
\right)
\to
\left(
\begin{array}{c}
\al^3 \\
\al^4 + \omega \al^1 \\
\al^1 \\
\al^2
\end{array}
\right).
$$
Decomposition (\ref{btDec}) is the $\si$-preimage
of (\ref{alDec}). Decomposition (\ref{siDec})
is computed by pairing $\dt_i$ with $\bt^k$ via the ordinary trace.
\end{proof}
\begin{corollary}
\label{decM}
The submodule of $\tau_{RE}$-symmetric tensors in $\M^{*\oo2}$ is
generated by the lowest weight vectors
$$\{\dt_1+\dt_2+\omega\dt_3, \dt_3+\dt_4-\omega\dt_2,\dt_5,\dt_6,\dt_7,\dt_8\}$$
The submodule of $\tau_{RE}$-antisymmetric tensors in $\M^{*\oo2}$ is
generated by the lowest weight vectors
$$\{\dt_1-\dt_2,\dt_3-\dt_4,\dt_9,\dt_{10}\}$$
\end{corollary}
\begin{proof}
Immediate from the last proposition.
\end{proof}
The permutation $\tau_{RE}$ can be reduced to the
submodule $\g^{\oo 2}\subset \M^{*\oo2}$. We study this problem in the
next section.
\section{Reducing $\tau_{RE}$ to $\g^{\oo 2}$}
The
left coadjoint module $\M^*$ contains a one-dimensional submodule
$\m_0$, which is  spanned by the invariant element $\sum_{i=1}^{n} q^{-2i+2}e^i_i$
and the $\g$-type submodule of $q$-traceless matrices (definition (\ref{PrTr}) 
of the quantum trace  will be given in the next section).
The decomposition  $\M^*=\m_0 \oplus \g$ leads to the decomposition
$$
\M^*\oo \M^*=\m_0 \oo \m_0 \;\oplus\; \m_0\oo \g \;\oplus\; \g \oo \m_0
\;\oplus\; \g\oo\g .
$$
The lowest weight vectors $\tilde\dt_3$ and  $\tilde \dt_4$  belonging to
$\g\oo\g$ are obtained by subtracting proper linear combinations of
$\dt_1$ and $\dt_2$  from $\dt_3$ and $\dt_4$:
\be{}
\label{tilde_si}
&\tilde \dt_3  =  \dt_3 - \frac{\omega}{1-q^{-2n}}\dt_1
- \frac{\omega}{1-q^{-2n}}\dt_2 , \quad
\tilde \dt_4  =  \dt_4 - \frac{\omega}{q^{2n}-1} \dt_1
- \frac{\omega}{1-q^{-2n}}\dt_2 .
\ee{}
Now we can evaluate the reduction $\tilde \tau_{RE}$ of
the permutation $\tau_{RE}$ to the submodule $\g^{\oo 2}$.
\begin{propn}
The composition
\be{}
\label{taug}
\g\oo\g  \to \M^*\oo\M^* \stackrel{\tau_{RE}} \longrightarrow
\M^*\oo\M^*  \to \g\oo \g,
\ee{}
where the left arrow is embedding and the right one
is the projection along $\m_0 \oo \M^* + \M^* \oo \m_0$,
defines an involutive $\U_q(\g)$-invariant permutation
$\tilde \tau_{RE}$ on $\g^{\oo 2}$.
The lowest weight vectors
\be{}
\tilde \dt_+  =  (1 - \frac{\omega q^{-n}}{n_q})\tilde \dt_3 +
(1 + \frac{\omega q^{ n}}{n_q})\tilde \dt_4,
&\quad &
\tilde \dt_-  = \tilde \dt_4 - \tilde \dt_3,
\label{pm}
\ee{}
where $\tilde \dt_i$ are introduced by (\ref{tilde_si}),
generate the $\g$-type submodules in $\g^{\oo 2}$ of
symmetric and antisymmetric tensors with respect to $\tilde \tau_{RE}$, correspondingly.
\end{propn}
\begin{proof}
The first statement is immediate.
It easy to check that the vector $\tilde \dt_+\in \g^{\oo2}\subset 
\M^{\oo2}$ is
$\tau_{RE}$-symmetric. As to the vector $\tilde \dt_-$,
it is the image of  $\tau_{RE}$-antisymmetric
vector $\dt_4-\dt_3$ under the projection $\M^{*\oo2}\to \g^{\oo2}$.
\end{proof}
\begin{corollary}
The submodules
of symmetric and antisymmetric tensors are generated
by the sets
$$
\{\dt_6-\frac{\omega}{1-q^{-2n}}\dt_5,\; \dt_7,\; \tilde\dt_+\},
\quad \mbox{and} \quad
\{\dt_9, \; \dt_{10},\; \tilde\dt_-\},
$$
of lowest weight vectors, respectively.
\end{corollary}
\begin{proof}
The vector $\dt_5$ is the square of the invariant element
$\sum_{i=1}^{n} q^{-2i+2}e^i_i$
spanning $\m_0$. The combination $\dt_6-\frac{\omega}{1-q^{-2n}}\dt_5$
is the Casimir element of $\g^{\oo 2}\subset \M^{*\oo 2}$.
The $\g$-type vectors $\tilde\dt_\pm$ were considered in the
above proposition.  Concerning the other lowest weight vectors,
they already lye in
$\g^{\oo2} \subset M^{*\oo2}$, so the proof is an immediate
consequence of  Corollary \ref{decM}.
\end{proof}
\begin{remark}
Unlike in the classical case, the submodule $\g^{\oo2}$ in $M^{*\oo2}$
is not preserved by the permutation  $\tau_{RE}$.  By this reason,
$\tilde\tau_{RE}$ is not a simple restriction of $\tau_{RE}$ to
$\g^{\oo2}$.
\end{remark}

As a left $\U_q(\g)$-module, the tensor square $\M^{*\oo 2}$
is endowed with another permutation coming from the universal
R-matrix representation:
$
A\oo B \to \Ru_2\tr B \oo \Ru_1\tr A, \quad A\oo B \in \M^{*\oo 2}.
$
It can be shown that its restriction to the $\g$-component
$\M^{*\oo2}_{\dt}\subset\M^{*\oo 2}$ acts according
to  the rule
$$
\left(
\begin{array}{c}
\dt_1 \\
\dt_2 \\
\tilde\dt_3 \\
\tilde\dt_4
\end{array}
\right)
\to
\left(
\begin{array}{r}
\dt_2 \\
\dt_1 \\
\tilde\dt_4 \\
q^{-2n}
\tilde\dt_3
\end{array}
\right).
$$
It is easy to see, that this operator has eigensubspaces
distinct from those of $\tau_{RE}$.
With this remark, we complete the study of the RE algebra. In
the remainder of the paper, we investigate relations between
$\A_{RE}({\M})$ and the quantization on $\g^*$ and its orbits.

\section{Algebra $\A_{RE}$ and the quantization on $\g^*$}
It was shown in \cite{Do1} that a two-parameter $\U_q(\g)$-covariant
quantization $\S_{q,t}(\g)$ of the Lie-Poisson bracket on $\g^*$
exists only in the $sl(n)$ case.
It is realized as the quotient algebra
of $\mbox{T}(\g)[t](q)$ by the quadratic-linear relations
\be{}
\label{quommutator}
x \oo y - \tilde\tau_{RE} (x \oo y) &=& t [x,y]_q, \quad
\ee{}
where $[\:\cdot, \cdot\:]_q\colon \g^{\oo2}\to \g$ is a deformed
Lie bracket. This is a $\U_q(\g)$-equivariant map $\g^{\oo2}\to \g$ sending
the submodule of $\tilde\tau_{RE}$-symmetric tensors to zero.
The commutator $[\:\cdot, \cdot\:]_q$ is uniquely defined up to a factor,
since the module of $\tilde\tau_{RE}$-antisymmetric tensors
contains the submodule isomorphic to $\g$ with multiplicity one.
In the limit $q\to 1$, the commutator turns into  the classical Lie
bracket. In this way, one recovers the  algebra $U(\g)[t]$ viewed as a quantization of the
Lie-Poisson bracket. Another limit $t\to 0$ yields a one-parameter family
$\S_{q,0}(\g)$, which can be interpreted as a
$\tilde \tau_{RE}$-commutative algebra of polynomials on $\g^*$.

The RE algebra is closely related to $\S_{q,t}(\g)$. To formulate that
relation, let us  introduce the quantum trace of the RE matrix $L$:
\be{}
\label{PrTr}
\Tr_q(L) &=& \sum_{i=1}^{n} q^{-2i+2} L^i_i =\Tr(DL),
\quad  \mbox{where}\quad
D =
\sum_{i=1}^{n} q^{-2i+2}e^i_i.
\ee{}
Remark that this definition makes sense for matrices with entries
being elements of any associative algebra $\A$ over $\C(q)$. The
quantum trace is the map
$\id\oo \Tr_q\colon \A\oo \M \to \A\oo\C(q) \sim \A$.  For the
generating matrix $L$ of the RE algebra, the elements $\Tr_q (L^k)$,
$k\in \N$, are invariant.  It follows from the commutation relations
(\ref{RER}) that they belong to the center of $\A_{RE}(\M)$.
\begin{theorem}
\label{SgAre}
The quotient algebra of $\A_{RE}(\M)$ by the ideal $(\Tr_q(L)-\la)$,
where $\la  \in \C$,
is isomorphic to the one-parameter sub-family
$t=\la \frac{ q-q^{-1}}{\mbox{\footnotesize Tr}_q(Id)}$
in $\S_{q,t}(\g)$.
\end{theorem}
\begin{proof}
In its essential part, this theorem was proven in \cite{Do1}. Here
we specify the exact relation between  the parameter $t$ and the value
of $\la$, using the explicit expression for the permutation $\tau_{RE}$
derived in Section \ref{qsl}, in terms of the lowest weight vectors.
By Proposition
\ref{tauRE}, the $\tau_{RE}$-symmetric $\g$-type submodules are
generated by the vectors $\dt_1-\dt_2$ and $\dt_4-\dt_3$. The first
one is irrelevant because it turns zero when projected to
$\g^{\oo2}\subset \M^{*\oo 2}$. The second vector is
represented as $\dt_4-\dt_3= \tilde\dt_4-\tilde\dt_3-\omega \dt_1$,
where $\tilde \dt_i$ belong to  $\g^{\oo2}$. Vector $\dt_1$ is
expressed through the matrix $D$ introduced in  (\ref{PrTr}) and the
lowest weight vector $e^n_1$ of the adjoint representation
in $\g$:  $\dt_1=
q^{-1}D\oo e^n_1$.  Comparing this with (\ref{quommutator}), we find
the value of the parameter $t$.
\end{proof}

\section{The quantization of coadjoint orbits}

It can be shown that  the
two-parameter family $\S_{q,t}(\g)$ can be restricted to any
semisimple orbit in $\g^*$. It this
connection, there arises the problem of explicit description of the
quantized manifolds in terms of ideals in $\S_{q,t}(\g)$. We solve
this problem (in a one-parameter setting) for certain classes of
orbits including the symmetric ones, using the relation between
$\A_{RE}$ and $\S_{q,t}(\g)$ (Theorem \ref{SgAre}).
Let ${\mathcal O}_A\subset \M$ be the orbit passing through the 
matrix $A$. Consider the two commutative diagrams 
$$
\begin{array}{ccc}
G\times \M &\longrightarrow &\M\\[8pt]
\uparrow &&\uparrow\\[8pt]
G\times \{A\}&\longrightarrow &{\mathcal O}_A
\end{array}
\quad\quad\quad
\begin{array}{ccc}
\mbox{Fun}(G)\oo \mbox{Fun}(\M) &\longleftarrow &\mbox{Fun}(\M)\\[8pt]
\id\oo \chi^A\downarrow  \quad\quad\quad \quad &&\downarrow\\[8pt]
\mbox{Fun}(G)\oo \C \;\;\;\quad\quad &\longleftarrow &\mbox{Fun}({\mathcal O}_A),
\end{array}
$$
where the horizontal arrows correspond to the action of the group $G$.
The right square represents morphisms of the polynomial algebras induced by
the maps of manifolds depicted on the left. The map $\chi^A$ is the 
character of the algebra $\mbox{Fun}(\M)$ corresponding to the point
$A\in \M$. We are going to quantize the right square; that will give
us realization of the quantized algebra  $\mbox{Fun}_q({\mathcal O}_A)$,
on the one hand, as a quotient of $\A_{RE}(\M)$ and as a subalgebra
in the quantized function algebra on $G$, on the other. 

Let $\F_q(G)$ be the Hopf algebra of quantized polynomial 
functions
on the group $G$. It is a quotient of the
FRT algebra by the additional relation $\det_q(T)=1$,
where $\det_q$ is the quantum determinant,
\cite{FRT}.
Since $\U_q(\g)$ acts on the RE algebra $\A_{RE}(\M)$, this action generates
the coaction of
$\F_q(G)$ on $\A_{RE}(\M)$. 
\begin{theorem}[\cite{KS}]
\label{KS}
Let $T^{-1}$ be the matrix with entries $\gm(T^i_j)$, where $T$ is
the generating matrix of the algebra $\F_q(G)$ and
$\gm$ the antipode.
The conjugation transformation $L\to T^{-1} L T$ of
an RE matrix $L$ with the FRT matrix $T$ whose entries commute
with the entries of $L$ is again an RE matrix.
\end{theorem}
\noindent 
It follows  that 
the correspondence  $L\to T^{-1} L T$ extends
to an homomorphism $\A_{RE}(\M)\to \F_q\oo \A_{RE}(\M)$.
Theorem \ref{KS} yields the quantization of the upper arrow
of the right square on the above diagram. To quantize the other 
maps, we replace $\chi^A$ by $\chi_q^A$, 
a character of the reflection
equation algebra. It is determined
by the correspondence $L^i_j  \to  A^i_j  \in \C(q)$
such that the numeric
matrix $A$ satisfies the 
reflection equation, \cite{KSS}, 
\be{} 
\label{RE1} A_2 S A_2 S &=& S A_2 S A_2 
\ee{}
supported in $\M^{\oo 2}$. 
Any solution to this equation gives rise to the 
algebra $\mbox{Fun}_q({\mathcal O}_A)$ closing
the commutative diagram
$$
\begin{array}{ccc}
\F_q(G)\oo \A_{RE}(\M)  &\longleftarrow &\A_{RE}(\M)\\[8pt]
\id\oo \chi^{A}_q\downarrow \;\quad\quad\quad\quad \quad &&\downarrow\\[8pt]
\F_q(G)\oo \C(q) \;\;\;\quad&\longleftarrow &\mbox{Fun}_q({\mathcal O}_A)
\end{array}
$$
so that the bottom arrow is embedding. Thus, we obtain 
\begin{theorem} \label{thmquan}
Let $A$ be a solution of the numeric RE (\ref{RE1}).
Then the algebra 
$\mbox{Fun}_q({\mathcal O}_A)$ in the diagram above is the quantization 
of the polynomial algebra on
the orbit ${\mathcal O}_A$ passing through the matrix $A$.
\end{theorem}

\begin{proof} It is clear that at $q=1$ the algebra $\mbox{Fun}_q({\mathcal O}_A)$ coincides
with the polynomial algebra on the orbit.
The flatness of $\mbox{Fun}_q({\mathcal O}_A)$ over $q$
follows from the fact that it is simultaneously a 
quotient and a subalebra of the flat $\C(q)$-algebras $\A_{RE}(\M)$
and $\F_q(G)$. 
\end{proof}

The matrix
$L(A)=T^{-1} A T$ possesses the following properties: 
\begin{lemma} For  any polynomial    function    ${\mathcal    P}$
in one variable,    
\be{}    \label{P}     
{\mathcal P}\bigl(L(A)\bigr)= L\bigl({\mathcal P}(A)\bigr). 
\ee{} The quantum trace is
invariant under  the  conjugation:  
\be{}  \label{Tr}  \Tr_q  (L(A))   =
\Tr_q(A).  
\ee{}  
\end{lemma}  

\begin{proof}  The  first  statement   is
evident. The second one is checked using the commutation  relations  in
the algebra $\F_q(G)$. 
\end{proof}

Note that solutions $A$ of (\ref{RE1})  yield  quantizations
which are quotients of one-parameter subfamilies in $\S_{q,t}(\g)$.
Those subfamilies are defined by paths 
in the parameter space $(q,t)$ crossing the axis $q=1$ at the point $t=0$. 
In particular, the solutions
with $\Tr_q(A)=0$ correspond to the path $t=0$. The limit $q\to 1$  for  $t$
separated from $0$ cannot be reached within our  approach.  
In the last subsections, we consider quantizations of various types 
of orbits, along the line of  Theorem \ref{thmquan}.

\subsection{Quantized symmetric orbits} 
To begin  with,
let us prove  a statement  relating  reflection  equation  algebras  in
different dimensions. Let $\M(n)$ and $\M(k)$ be the matrix  algebra  of
$n\times n$ and $k\times k$ matrices, $0 < k <n$.  The  homomorphism
$\M(k)\to \M(n)$  of  embedding  as  the  left  upper  block   induces   a
contravariant epimorphism of the function algebras. The  same  holds  in  the
quantum situation: 
\begin{propn} 
\label{reduction} 
The quotient  algebra of 
$\A_{RE}\bigl(\M(n)\bigr)$ by the relations $L^i_j=0, \quad i>k  \;\mbox{or}
\;  j>k,$  is  isomorphic  to  $\A_{RE}\bigl(\M(k)\bigr)$.   
\end{propn}
\begin{proof} Let $P^+$  and  $P^-$  be  the  projectors  from  $\C^n  =
\C^k\oplus \C^{n-k}$ to the first and the second addends,  respectively.
Denote $R^{++}$ the projection of the R-matrix  (\ref{R})  to  $\M(k)\oo
\M(k)$. Up to a nonzero scalar factor, this is the R-matrix in
dimension $k$. The equalities 
\be{}
(P^+\oo 1)  R = R^{++} + q^{-\frac{1}{ n}} P^+\oo P^-, 
\label{+R}
\\
R  (1\oo P^+) = R^{++} + q^{-\frac{1}{n}} P^- \oo P^+
\label{R+}
\ee{}
follows directly from (\ref{R}). 
The matrix $L$ is equal to $L^+=P^+LP^+$ modulo the ideal  specified  in
the hypothesis. In terms of the matrix $R$ instead of the  Hecke  matrix
$S$,   the reflection equation is rewritten 
as 
$$
R_{21} L^+_1 R L^+_2 = L^+_2 R_{21} L^+_1 R
$$
modulo the relations of concern.
Using (\ref{+R}) and (\ref{R+}), the matrix $R$ can be replaced by 
$R^{++}$ thus leading to the RE in dimension $k$.
\end{proof} 
Proposition \ref{reduction} suggests a method of  building
solutions to (\ref{RE1})  by those  in  smaller  dimensions
extending them by zero matrix elements. 
For example,    the     projectors
$P_k=\sum_{i=1}^{k}e^i_i$ of rank $k=1,\ldots, n$  satisfy  (\ref{RE1}).
Indeed, each can be obtained from the $k\times k$ unit matrix (which  is
apparently a solution to the RE) by extending it with  zeros  to  the  $n\times  n$
matrix. 

Let  us  introduce  the  quantum
integer numbers 
\be{} 
\hat k_q =  \sum_{i=1}^{k}  q^{-2i+2}= \Tr_q(P_k).  
\ee{}
\begin{theorem}  \label{qGr}
The quotient algebra of $\A_{RE}(\M)$ by the  relations  
\be{}  \label{gr}
L^2=L,\quad  Tr_q(L)=\hat  k_q  
\ee{}  
yields   a   $\U_q(\g)$-covariant
quantization $\mbox{Fun}_q({\mathcal O}_{P_k})$ 
of the symmetric orbit ${\mathcal O}_{P_k}$ passing through
the projector $P_k$
\end{theorem}  
\begin{proof}  Relations  (\ref{gr})
follow from identities (\ref{P}) and (\ref{Tr}) as applied to $A=P_k$. 
The projectors
$P_k$  are  stable  under  the  right  adjoint  action   (\ref{ad})   of
$\U_q(sl(k))\oo \U_q(sl(n-k))\oo \U_q(sl(1))$ as a quantum subgroup in  $\U_q(sl(n))$.
Therefore, the subalgebra  in  $\F_q(G)$  generated  by  the  RE  matrix
$L(P_k)$ is invariant under the right action  of 
$\U_q(sl(k))\oo \U_q(sl(n-k))\oo \U_q(sl(1))$. 
\end{proof} 
\begin{remark}  
The  quantum ${\mathcal O}_{P_1}$ may be obtained directly from  the  
description
of $\A_{RE}$ given in Section \ref{SecRE}.  Indeed,  imposing
the additional conditions $V^{\wedge 2}=0$ and $V^{*\wedge 2}  =  0$  one
comes to the subalgebra in $\Sym_q(V)\check\oo \Sym_q(V^*)$ generated by
$e_i\oo f^j$ with $e_i$ and $f^j$ commutative in the sense of permutations
(\ref{VV}--\ref{VV*}).  The matrix elements $L^j_i=e_i f^j $ satisfy 
the equality $ L^l_i L^j_l = e_i  (f^l e_l) f^j = Tr_q(L) L^i_j  $  
following  from  the  commutation
relations.  So,  conditions  (\ref{gr})  for   the   case   $k=1$   hold
simultaneously  in  this  algebra. Another description of the quantum projective
space is found in \cite{DGK}.  
\end{remark}  
As  in  the  classical
situation, we can consider solution $\la P_k$  with  arbitrary 
$\la\not = 0$. This will lead to the  relation  $L^2=\la  L$  and  the
corresponding rescaling of the quantum trace. The resulting conditions give  the
deformation quantization of the orbit passing through $\la P_k$.  It  may
be regarded as the quantization of the  same  manifold  but along the
path  $  t=\la\omega \frac{\hat  k_q}{\hat  n_q}$   in   the
two parameter space of  the  universal  family  $\S_{q,t}(\g)$  (Theorem
\ref{SgAre}). 

\subsection{The Cayley-Hamilton identity} The goal of this
subsection is to exhibit correspondence between 
Theorem  \ref{qGr} with identities in  $\A_{RE}(\M)$
of the Cayley-Hamiltonian type, \cite{PS}. In the algebra  $\A_{RE}(\M)$,
the $k$-th powers $L^k$ of the generating matrix form the left coadjoint
module $\M^*$.  The  quantum  traces  $\Tr_q(L^k)$,  $k=1,2,\ldots$  are
$\U_q(sl(n))$-invariant and central, (see, e.g., \cite{PS}).  
\begin{theorem}[\cite{PS}]  The  generating  matrix   of   the
reflection  equation  algebra  $\A_{RE}(\M)$  obeys  the   relation   
$$
\sum_{j=0}^n \si^j_q(L)(-L)^{n-j} =0, 
$$ where $\si^j_q(L)$ are  central
elements expressed through the quantum traces of powers in  $L$  by  the
recursive  formula  
\be{}  \hat  k_q  \si^k_q(L)  &=&   \sum_{j=0}^{k-1}
\si^j_q(L)(-1)^{k-j+1}\Tr_q(L^{k-j}), \quad \si^0_q(L)=1.  
\label{sym_f}
\ee{}  
In  particular,  $\si^{n+1}_q(L)=0$.  
\end{theorem}  
It follows  that  imposing  the  projector  condition  on  the  matrix  $L$
specifies  the  quantum  trace  modulo  the   finite   set   of   values.
\begin{corollary} 
The condition $L^2=L$ on the generating matrix of  the
$\A_{RE}(\M)$ implies the discrete set $\{\hat 0_q,\hat 1_q,\ldots, \hat
n_q\}$ of values for the quantum  trace.  
\label{Trace}  
\end{corollary}
\begin{proof} The  statement  follows  from  the  formula  (\ref{sym_f})
since,  under  the  hypothesis   made,   one   has   $$   \si^{n+1}_q(L)=
\frac{\Tr_q(L)\bigl(\Tr_q(L)-\hat  1_q\bigr)\ldots   \bigl(\Tr_q(L)-\hat
n_q\bigr)}{\widehat{(n+1)_q}!}=0. $$ 
\end{proof} 
All the possible values
of $\Tr_q(L)$ are realized, by Theorem  \ref{qGr},
giving flat deformations of symmetric spaces.  
Note that the  Cayley-Hamiltonian identity in the RE algebra
was directly used for construction of the quantum sphere ${\mathbb S}^2_q$ 
in \cite{GS}.

\subsection{On the quantization of  non-semisimple  orbits}  
As was already mentioned, the
two-parameter quantization $\S_{q,t}(sl(n))$ may  be  ``restricted''  to
every semisimple orbits. There are no definite assertions of  that  kind
concerning orbits passing through nilpotent elements. In  this  section,
we prove 
\begin{theorem} There exists a  one-parameter  $\U_q(\g)$-covariant
quantization of nilpotent orbits in $\M$  satisfying  the matrix
equation $A^2=0$. It is a restriction of the one-parameter subfamily
 $t=0$  in  the  universal
two-parameter quantization $\S_{q,t}(\g)$. 
\end{theorem} 

\begin{proof} Like  in
the case of symmetric orbits, we seek  for  a  solution  to  the  numeric
reflection equation, which will  realize  the  quantized  algebra  as  a
subalgebra in  $\F_q(G)$.  Specializing  to  the  skew-diagonal  matrices,
$$
A=\sum_{i=1}^n\la_i e^{i'}_i,  \quad  i'=  n+1-i,
$$  
we  find it in the form $\la_i \la_{i'}=\la^2, \quad i=1,\ldots,n.$  In
case of non-zero $\la$  we  come  to  nondegenerate  matrices  with  two
eigenvalues $\pm\la$. These solutions were  found  in  \cite{KSS}.  They  lead  to
other realization  of  quantum  symmetric  orbits  than  by  means  of
projectors. In the case of $\la=0$,
these matrices being squared  are  zero.  By  re-enumerating  the  basis
elements, they can be brought to the sum of $2\times 2$ jordanian blocks
yielding all such matrices. The last statement of the theorem
holds because $\Tr_q(A)=0$.
\end{proof}

\subsection{On the quantization of bisymmetric orbits} To quantize 
the symmetric orbits in $g^*$, we used the projectors, i.e.,  
the  semisimple
elements  with  the eigenvalues  $1$  and  $0$.  There are non-degenerate
semisimple solutions to the matrix reflection equation with two eigenvalues. 
For example, one can take $\sum_{i=1}^n e^i_{n+1-i}$. Such matrices can  be
used for  constructing  solutions with  three  eigenvalues  by the  embedding
method, along the line of Proposition \ref{reduction}. In this way,  one
gains the additional zero eigenvalue. One might have expected that this
will provide a tool for quantizing all the semisimple  orbits,  which  are
classified as homogeneous spaces by  the  number  of  eigenvalues  and  
their  multiplicities.
However, that  is probably impossible, due to certain indications.  Indeed,  
in \cite{KSS}, there were written out all  non-degenerate  solutions  up  to
$n=4$. They have at most two eigenvalues and lead  to  different
quantizations for the symmetric orbits than by means of projectors. Among 
them, there
are also  q-traceless  matrices.  We  conjecture  that non-degenerate
matrix solutions  to the  reflection  equation with arbitrary
value of the quantum trace do  exist  for
each symmetric orbit. This implies the existence of the quantization 
via the RE algebra 
characters for  every
bisymmetric orbit, i.e., consisting of matrices with three eigenvalues.
Let us present,  following  \cite{KSS},  the  non-degenerate  q-traceless
RE matrices up to $n=4$. 
$$ A^{1,1}=
\left( 
\begin{array}{cc}
   & 1  \\
 1 & 
\end{array} 
\right) ,
\quad A^{2,1}=
\left(  
\begin{array}{ccc}
-q^{-2} & & \hat 2_q \\ 
& 1 &  \\  
1  &  &  
\end{array}  
\right)  
$$  
$$
A^{2,2}= 
\left( 
\begin{array}{cccc}
 & & & 1 \\
 & &1& \\ 
 &1 & & \\ 
1 &   &   &    
\end{array}   
\right)    ,\quad    
A^{3,1}= 
\left(
\begin{array}{cccc}
-q^{-2}\hat 2_q& & & \hat 3_q \\
 & 1 & & \\
 & & 1 & \\
 1 & & &
\end{array}
 \right). 
$$ 
The matrix $A^{2,2}$ is interesting  from  the
physical point of view because it yields the quantization of the
twistor space.
Using these RE matrices together with the projectors $P_k$, we  can  
cover all the semisimple orbits for $n\leq 4$, excepting the  maximal 
orbit in dimension
$n=4$. But the maximal orbits can  be  quantized  by  specifying  values  
of the Casimir elements, i.e., $\Tr(L^i)$, $i=1,\ldots,n-1$.  So  we 
conclude that, at least to dimension  $n=4$, all the semisimple orbits 
can  be  quantized  explicitly 
as quotients of the reflection equation algebra.

\small

\bibliographystyle{amsalpha}

\begin{thebibliography}{A}

\bibitem[AFS]{AFS} A. Alekseev, L. Faddeev, M. Semenov-Tian-Shansky,
\textit{ Hidden quantum group inside Kac-Moody algebra,}
     Proceedings of the  Euler  International  Mathematical
   Institute on Quantum Groups, 
   Lect. Notes Math.  {\bf  1510}  (Springer, Berlin, 1992) 148. 
\bibitem[Cher]{Cher}  I. V. Cherednik,
\textit{ Factorizing particles on a half line, and root systems,}  
Teoret. Mat. Fiz. {\bf 61} (1984), \# 1, 35--44.
\bibitem[D]{D} V.  G.  Drinfeld,  
  \textit{  Quantum Groups}, in Proc. Int. Congress of Mathematicians, 
  Berkeley,  1986,  ed. A.V. Gleason, AMS, Providence (1987) 798. 
\bibitem[Do]{Do} J. Donin,
\textit{ Double quantization on the coadjoint representation of $sl(n)$,}
Czech J. of Physics, {\bf 47}, n.11, 1997, 1115-1122.
\bibitem[Do1]{Do1}  J. Donin, 
   \textit{ $U_h(\g)$-invariant  quantization  of  
   coadjoint  orbits  and vector bundles over  them,}  
   Preprint  54  (2000),  Max-Plank-Institute,
J. of Geometry and Physics, {\bf 38} (2001) 54.
\bibitem[DoG]{DoG}  J. Donin,  D. Gurevich,
  \textit{ Some Poisson structures associated to Drinfeld-Jimbo 
  R-matrices and their quantization,}
 Israel J. of Math. {\bf 92} (1995) 23.
\bibitem[DGK]{DGK} J. Donin, D. Gurevich, S. Khoroshkin, 
 \textit{ Double quantization of $\C P^n$ type by  generalized
Verma   modules,}   Preprint   math.QA/9803155,
J. of Geom. and Phys., v.28, 1998, 384.   
\bibitem[FRT]{FRT} L. Faddeev,  N. Reshetikhin, and L. Takhtajan, 
 { \em Quantization of Lie groups and Lie algebras}  
 Leningrad Math. J. {\bf 1}  (1990)  193.  
\bibitem[G]{G} D. Gurevich, 
  \textit{   Algebraic aspects of quantum Yang-Baxter equation,}  
  Leningrad Math. J. {\bf 2} (1991) 802.
\bibitem[GS]{GS} D. Gurevich, P. A. Saponov \textit{  Quantum  sphere  via  reflection  equation
algebra,} Preprint math.QA/9911141. 
\bibitem[KT]{KT}  S. M. Khoroshkin, V. N. Tolstoy 
\textit{  Universal R-matrix for quantized (super) algebras,}
             Comm. Math. Phys. {\bf 141} (1991) 559. 
\bibitem[KR]{KR} A. N. Kirillov, N. Yu. Reshetikhin,
         \textit{ q-Weyl group an multiplicative formula for
          R-matrices,}
            Comm. Math. Phys. {\bf 130} (1990) 421. 
\bibitem[K]{K} P. P. Kulish, 
  \textit{ Quantum groups, q-oscillators, and covariant algebras,} 
  Theor. Math. Phys. {\bf 94} (1993) 193.
\bibitem[KSkl]{KSkl} P. P. Kulish, E. K. Sklyanin 
  \textit{Algebraic  structure  related to the reflection equation},
  J. Phys. A {\bf  25} (1992)  5963.  
\bibitem[KS]{KS}  P.  P.  Kulish,  R.  Sasaki   
  \textit{ Covariance properties of reflection  equation  algebras,}  
  Prog.  Theor. Phys. {\bf 89} $\# 3$  (1993)  741.  
\bibitem[KSS]{KSS}  P. P.  Kulish,  R. Sasaki, and C. Schweibert,
  \textit{ Constant solutions  of  reflection  equations  and quantum groups,}
                 J.  Math.  Phys  {\bf  34}  $\#  1$  (1993)  286.
\bibitem[Lu]{Lu} G. Lusztig, ``Introduction to quantum groups'', 
Progress in Mathematics, 110. Birkh.
Boston, Inc., Boston, MA, 1993. 
\bibitem[Mj]{Mj} S. Majid, ``Foundations of  quantum
group theory'', Cambridge University  Press,  1995.  
\bibitem[PS]{PS}  P.  N. Pyatov,  P.  A.  Saponov  
  \textit{  Characteristic  relations  for   quantum matrices,} 
  J. Phys. A {\bf 28} $\# 3$ (1995) 4415. 
\bibitem[Ros]{Ros} M. Rosso,
         \textit{ An analog of P.B.W. theorem and the universal
        $U_h sl(N+1)$,}
        Comm. Math. Phys. {\bf 124} (1991) 307. 
 
\end{thebibliography}

\end{document}